\documentclass[
12pt,
authoryear,
times
]{elsarticle}

\usepackage{graphicx,color}

\usepackage[colorlinks=true, breaklinks=true, pdfstartview=FitV, linkcolor=red, citecolor=blue, urlcolor=black]
{hyperref}

\usepackage[T1]{fontenc}

\usepackage{algorithm}
\usepackage{algorithmicx}
\usepackage{algpseudocode}

\newcommand{\Vec}[1]{\mbox{\boldmath$#1$}}
\journal{Celestial Mechanics and Dynamical Astronomy}

\begin{document}

\begin{frontmatter}

\title{ Note on the ideal frame formulation 
 }

\author[roa]
{ Martin Lara\fnref{SDG} 
}
\ead{mlara0@gmail.com}


\fntext[SDG]{ GRUCACI, University of La Rioja, and Space Dynamics Group -- UPM  }

\address[roa] {C/ Luis de Ulloa, s.n., 26004 Logro\~no, Spain}

\date{}
\begin{abstract}
An implementation of the ideal frame formulation of perturbed Keplerian motion is presented which only requires the integration of a differential system of dimension 7, contrary to the 8 variables traditionally integrated with this approach. The new formulation is based on the integration of a scaled version of the Eulerian set of redundant parameters, and slightly improves runtime performance with respect to the 8 dimensional case while retaining comparable accuracy.
\end{abstract}


\end{frontmatter}

\section{Introduction}

In the times in which computers were humans and the success in ephemeris computation was dependent on the available resources for hiring computers, astronomers made serious efforts in developing efficient formalisms and methods for the numerical integration of orbits. Since then, the advantages of Encke's or Hansen's formulations over the straightforward Cowell's approach were mandatorily described in classic textbooks on celestial mechanics \cite[see chap.~11 of][for instance]{Danby1992}. However, in an epoch of computational plenty, one may think that ``computing technology has overcome the need to use such techniques'' \cite[p.~516]{Vallado2007}, thus relegating them into oblivion. Quite on the contrary, current needs for orbit propagation, as for instance, those related to Space Situational Awareness, motivate renewed interest in the efficient numerical integration of orbital motion, a fundamental part of which is related to finding the more suitable formulation of the equations of motion.
\par

The identification of slow and fast effects of a force can be beneficial in numerical integration. Indeed, it helps in devising formalisms that increase the speed of the numerical integration by abbreviating the computations, on the one hand, and keep as much significant digits as possible by computing the variations of the slow varying elements with respect to a suitable inertial frame, on the other. For its part, the accuracy of the fast variables will also be increased because it generally depends on the precision achieved in the determination of the slow variables. 
\par

Different approaches to the separation of slow and fast effects are widely enclosed under the name of \emph{variation of parameters} methods (VOP) ---the ``parameters'' being particular combinations of the integration constants of the Keplerian motion--- each of which has merits and drawbacks. Among them, those based on the fundamental role played by the orbital plane when the component of the external forces in its normal direction is small, offer clear advantages. In that case, the orbital motion can be viewed as the composition of two different effects: the slow rotation of the orbital plane and the fast rotation of the particle within that plane. Even though the differential system representing the motion is not decoupled, one can come very close to the separation when using Hansen's \emph{ideal frame} concept \citep{Hansen1857}.
\par

The ideal frame is a moving frame attached to the orbital plane with the remarkable characteristic that the velocity is the same when measured in the ideal frame as when measured in the inertial frame. Besides, the attitude of the ideal frame with respect to the fixed frame can be materialized by the set of Eulerian redundant parameters, in this way avoiding singularities related to the evaluation of circular functions \citep{Musen1958}. Additional benefits of the ideal frame formulation are obtained when the motion in the orbital plane is decomposed into the slow evolution of the ellipse as given by the attitude of the apsidal frame, on the one hand, and the timing on the osculating ellipse, on the other \citep{Deprit1975}. The latter is as well measured in a slow scale after standard regularization \citep{SharafAwadNajmuldeen1992,DepritElipeFerrer1994,PalaciosCalvo1996}. Approaches based on Hansen's ideal frame give rise to extremely efficient formulations of the orbital motion, in which the use of Deprit's  \emph{ideal elements}, given by the projections of the eccentricity vector in the ideal frame, shows specially useful. An updated list of references using these formalisms can be found in \citep{UrrutxuaSanjurjoRivoPelaez2016}. 
\par

In addition to speed and significance, the non-singular character of the variables and the reduced dimension of the differential system to integrate are normally listed as desirable characteristics of a given formulation \citep{Fukushima2007}. All of them can be achieved within the ideal frame approach, which, besides, in the finer formulations, only needs to integrate 8 redundant variables. Still, it will be shown that the ideal frame formulation can be reduced to the integration of just 7 variables without loosing any of its recognized merits. The modification of the standard ideal frame algorithm is very simple. It consists in replacing the set of Eulerian parameters by a modified set in which each Eulerian parameter is scaled by the square root of the modulus of the angular momentum vector, thus avoiding the need of integrating the variation of this scalar. 
\par

For completeness, the ideal frame formulation is recalled in Section \ref{s:iff}, closely adhering to the elegant account in \citep{Deprit1975}. The particular modifications of the formalism that lead to the dimension reduction of the differential system are introduced in Section \ref{ss:modEulerp}, and the regularized systems are formulated in the new variables in following sections. Finally, Section \ref{s:performance} summarizes some test results which show that the new formulation may lead to a slightly better performance in terms of runtime when compared to the traditional formulation in 8 dimensions.
\par

\section{Ideal frame formulations} \label{s:iff}

Let $O$ be a fixed point in space, and let $\mathcal{S}\equiv\mathcal{S}(O,\Vec{i},\Vec{j},\Vec{k})$, the \emph{space frame}, be an inertial frame in which $\Vec{x}$ denotes the position with respect to $O$ of a mass particle that is moving with velocity $\Vec{X}$. In a Newtonian field, the equations of motion are written in the form of the first order differential system
\begin{eqnarray} \label{xxp}
\dot{\Vec{x}} &=& \Vec{X}, \\  \label{xpXp}
\dot{\Vec{X}} &=& \Vec{F}\equiv-\frac{\mathcal{G}M}{r^3}\Vec{x}+\Vec{P},
\end{eqnarray}
where $\Vec{F}$ represents the force (per unit of mass), $\mathcal{G}M$ is the gravity constant of the central body, $r=\|\Vec{x}\|$, and the perturbing force $\Vec{P}\equiv\Vec{P}(\Vec{x},\Vec{X},t)$ will be assumed to be small when compared to the Keplerian attraction.
\par

The time evolution of the perturbed Keplerian motion can be achieved by direct integration of Eq.~(\ref{xpXp}) in the space frame from proper initial conditions $\Vec{x}_0=\Vec{x}(t_0)$, $\Vec{X}_0=\Vec{X}(t_0)$, a straightforward approach that is customarily known as \emph{Cowell's formulation}. However, in this formalism all variables may vary fast and, for that reason, numerical integration methods need smaller step sizes with the consequent increase of the computational load and faster accumulation of truncation errors. The classical alternative is to integrate the differences between the perturbed orbit and some reference Keplerian orbit, the so-called Encke's formulation. As far as these differences remain small, the numerical integration can proceed with considerably larger step sizes than in Cowell's formulation. Advantages and inconveniences of both formalisms are profusely described in classical textbooks on celestial mechanics.
\par

Let
\begin{equation} \label{xxX}
\Vec{G}=\Vec{x}\times\Vec{X},
\end{equation}
be the (instantaneous) angular momentum vector per unit of mass, which defines the instantaneous \emph{orbital plane}, and let
\begin{equation} \label{Kepler:ev}
\Vec{e}=\frac{1}{\mathcal{G}M}\Vec{X}\times\Vec{G} -\frac{1}{r}\Vec{x},
\end{equation}
be the (instantaneous) eccentricity vector. The unit vectors $\Vec{n}=\Vec{G}/G$, $\Vec{a}=\Vec{e}/e$, and $\Vec{b}=\Vec{n}\times\Vec{a}$, where $G=\|\Vec{G}\|$ and $e=\|\Vec{e}\|$, define the \emph{apsidal frame} $\mathcal{A}\equiv\mathcal{A}(O,\Vec{a},\Vec{b},\Vec{n})$. Alternatively to the integration of Eqs.~(\ref{xxp}) and (\ref{xpXp}), the time evolution of the perturbed Keplerian motion is described by the motion of $\mathcal{A}$ with respect to $\mathcal{S}$ together with the evolution of $\Vec{x}$ with respect to the apsidal frame. That is, the integration of $\dot{\Vec{G}}$ and $\dot{\Vec{e}}$, the time variations of the instantaneous orbital plane and the instantaneous ellipse in the orbital plane, respectively, plus an additional equation for the timing \citep[see, for instance,][and references therein]{RoyMoran1973}. The dimension of the differential system has been raised from six to seven, but the constraint $\Vec{G}\cdot\Vec{e}=0$ reduces the flow to the correct dimension and can be used to check the quality of the numerical integration.
\par

This alternative approach takes benefit of the slow variation of the elements $\Vec{G}$ and $\Vec{e}$, which are constants of the (unperturbed) Keplerian motion, and is commonly known as the \emph{variations of parameters method} or VOP. Different parameters can be used, and there is no need of rising the dimension of the differential system \citep[see chap.~10 of][for instance]{Battin1999}.
\par

\subsection{The orbital frame}

The \emph{orbital frame} $\mathcal{O}\equiv\mathcal{O}(O,\Vec{u,v,n})$ is defined by the unit vectors: $\Vec{u}$ in the particle's direction, $\Vec{n}$ in the direction of $\Vec{G}$, and $\Vec{v}$ completing a direct orthonormal frame, viz.
\begin{equation} \label{orbitalframe}
\Vec{u}=\frac{\Vec{x}}{r}, \qquad \Vec{n}=\frac{\Vec{G}}{G}, \qquad \Vec{v}=\Vec{n}\times\Vec{u}.
\end{equation}
The space and orbital frames are linked by a rotation defined by the matrix
\begin{equation} \label{SO}
\mathcal{M}\equiv 
\left(\begin{array}{ccc} 
\Vec{u}\cdot\Vec{i} & \Vec{v}\cdot\Vec{i} & \Vec{n}\cdot\Vec{i} \\
\Vec{u}\cdot\Vec{j} & \Vec{v}\cdot\Vec{j} & \Vec{n}\cdot\Vec{j} \\
\Vec{u}\cdot\Vec{k} & \Vec{v}\cdot\Vec{k} & \Vec{n}\cdot\Vec{k} 
\end{array}\right).
\end{equation}
That is, for a given vector $\Vec{q}$ expressed in the orbital frame as the column matrix $\Vec{q}_\mathcal{O}$, corresponding coordinates in the space frame are computed as $\Vec{q}_\mathcal{S}=\mathcal{M}\Vec{q}_\mathcal{O}$.
\par

The motion of the orbital frame with respect to a fixed frame is given by the time variation of Eq.~(\ref{orbitalframe}), which is computed as follows. First of all we recall that, because the derivative of a unit vector is orthogonal to itself, if $\Vec{q}=q\Vec{w}$, where $\|\Vec{w}\|=1$, it can be checked that \citep{DepritElipeFerrer1994}
\begin{equation} \label{qwp}
\dot{q}=\dot{\Vec{q}}\cdot\Vec{w}, \qquad
\dot{\Vec{w}}=q^{-2}(\Vec{q}\times\dot{\Vec{q}})\times\Vec{w}.
\end{equation}
Applying this rule to $\Vec{x}=r\Vec{u}$ and to $\Vec{G}=G\Vec{n}$, and in view of $\dot{\Vec{G}}=\Vec{x}\times\Vec{F}$, as derived from Eq.~(\ref{xxX}), one easily finds
\begin{eqnarray} \label{rp}
\dot{r} &=& \Vec{X}\cdot\Vec{u}, \\
\label{Gp}
\dot{G} &=& r(\Vec{F}\cdot\Vec{v})=r(\Vec{P}\cdot\Vec{v}),
\end{eqnarray}
and
\begin{eqnarray} \label{vup}
\dot{\Vec{u}} &=& \frac{G}{r^2}\Vec{v}, \\ \label{vnp}
\dot{\Vec{n}} &=&-\frac{r}{G}(\Vec{P}\cdot\Vec{n})\Vec{v}, \\ \label{vvp}
\dot{\Vec{v}} &=&-\frac{G}{r^2}\Vec{u}+\frac{r}{G}(\Vec{P}\cdot\Vec{n})\Vec{n},
\end{eqnarray}
where the latter has been obtained after differentiation of $\Vec{v}$ in Eq.~(\ref{orbitalframe}).
\par

Integration of Eqs.~(\ref{Gp})--(\ref{vvp}) require to know the time variation of $r$. It is obtained from differentiation of Eq.~(\ref{rp}) taking into account that $\Vec{X}=\dot{r}\Vec{u}+r\dot{\Vec{u}}$, and the following use of Eqs.~(\ref{xpXp}) and (\ref{vup}), viz.
\begin{equation} \label{ddr}
\ddot{r} =\frac{G^2}{r^3}+\Vec{F}\cdot\Vec{u} =\left(\frac{G}{r}-\frac{G}{p}\right)\frac{G}{r^2}+\Vec{P}\cdot\Vec{u},
\end{equation}
where $p=G^2/(\mathcal{G}M)$ is the orbit parameter.
\par

The differential system made of Eqs.~(\ref{Gp})--(\ref{ddr}) is of dimension twelve, but it accepts the six constraints
\begin{equation} \label{constraints}
\|\Vec{u}\|=\|\Vec{v}\|=\|\Vec{n}\|=1, \qquad 
\Vec{u}\cdot\Vec{v}=\Vec{v}\cdot\Vec{n}=\Vec{n}\cdot\Vec{u}=0,
\end{equation}
derived from the orthogonality conditions, and hence it remains equivalent to Eqs.~(\ref{xxp})--(\ref{xpXp}).
\par

\subsection{Ideal frames}

The motion of the orbital frame defined by Eqs.~(\ref{vup})--(\ref{vvp}) can be viewed as the rotation
\begin{equation} \label{upvpnp}
\dot{\Vec{u}}=\Vec\omega\times\Vec{u}, \qquad
\dot{\Vec{v}}=\Vec\omega\times\Vec{v}, \qquad
\dot{\Vec{n}}=\Vec\omega\times\Vec{n}, 
\end{equation}
with the angular velocity
\begin{equation} \label{omega}
\Vec\omega=\frac{r}{G}(\Vec{P}\cdot\Vec{n})\Vec{u}+\frac{G}{r^2}\Vec{n}.
\end{equation}
That is, the composition of a rotation with the angular velocity
\begin{equation} \label{wI}
\Vec{\omega}^*=\frac{r}{G}(\Vec{P}\cdot\Vec{n})\Vec{u},
\end{equation}
from the fixed frame $\mathcal{S}$ to an intermediate frame $\mathcal{I}\equiv\mathcal{I}(O,\Vec{u}^*,\Vec{v}^*,\Vec{n})$, followed by a rotation from $\mathcal{I}$ to the orbital frame $\mathcal{O}$ with the angular velocity
\begin{equation} \label{wK}
\Vec{\omega}_K=\frac{G}{r^2}\Vec{n}.
\end{equation}

Because the position of $\mathcal{O}$ with respect to $\mathcal{I}$ is given by the Keplerian rotation $\Vec\omega=\Vec{\omega}_K$, the system $\mathcal{I}$ enjoys a remarkably property, which is easily derived from the theorem of the moving frame:
\begin{equation} \label{darboux}
\dot{\Vec{x}}=\frac{\partial\Vec{x}}{\partial{t}}+\Vec\omega^*\times\Vec{x},
\end{equation}
where the partial derivative notation is used with the meaning of differentiation in the rotating frame. Then, because $\Vec\omega^*\times\Vec{x}=\Vec{0}$ it happens that the velocity is the same when measured in $\mathcal{I}$ or in the inertial frame. Because of that, the rotating frame $\mathcal{I}$ is customarily called \emph{ideal}.
\par

In a perturbation problem $\Vec{P}$ is small when compared to the Keplerian attraction and so it is $\Vec{\omega}^*$ when compared to $\Vec{\omega}_K$. Hence, the motion of the particle can be viewed as a slow rotation of the orbital plane followed by the fast motion of the particle in this plane, the latter being a planar problem. This fact makes quite relevant the study of the motion of the ideal frame
\begin{eqnarray} \label{upI}
\dot{\Vec{u}}^* &=& \Vec{\omega}^*\times\Vec{u}^*, \\ \label{vpI}
\dot{\Vec{v}}^* &=& \Vec{\omega}^*\times\Vec{v}^*, \\ \label{npI}
\dot{\Vec{n}}\phantom{^*} &=& \Vec{\omega}^*\times\Vec{n}.
\end{eqnarray}
Furthermore, to avoid loss of significant digits in the integration of Eqs.~(\ref{upI})--(\ref{npI}), it is advisable to refer the integration to such a fixed frame that $\mathcal{I}$ remains as close to it as possible. Thus, the \emph{departure frame} $\mathcal{D}\equiv\mathcal{D}(O,\Vec{u}_0,\Vec{v}_0,\Vec{n}_0)=\mathcal{O}(O,\Vec{u}(t_0),\Vec{v}(t_0),\Vec{n}(t_0))$ is customarily chosen as the fixed frame, which is linked with $\mathcal{S}$ by a constant rotation defined by the matrix $\mathcal{M}_0= \mathcal{M}(t_0)$ computed from Eq.~(\ref{SO}).
\par

Let $\theta^*$ be the angle from $\Vec{u}^*$ to the position vector $\Vec{x}=r\Vec{u}$, so that
\begin{eqnarray} \label{uI}
\Vec{u} &=& \Vec{u}^*\cos\theta^*+\Vec{v}^*\sin\theta^*, \\ \label{vI}
\Vec{v} &=& \Vec{v}^*\cos\theta^*-\Vec{u}^*\sin\theta^*.
\end{eqnarray}
Then,
\begin{equation} \label{wIss}
\Vec{\omega}^*=\frac{r}{G}(\Vec{P}\cdot\Vec{n})(\Vec{u}^*\cos\theta^*+\Vec{v}^*\sin\theta^*),
\end{equation}
and
\begin{eqnarray} \label{duI}
\dot{\Vec{u}}^*  &=& -\frac{r}{G}(\Vec{P}\cdot\Vec{n})\Vec{n}\sin\theta^*, \\ \label{dvI}
\dot{\Vec{v}}^* &=& \frac{r}{G}(\Vec{P}\cdot\Vec{n})\Vec{n}\cos\theta^*, \\ \label{dnI}
\dot{\Vec{n}}\phantom{^*} &=& -\frac{r}{G}(\Vec{P}\cdot\Vec{n})\Vec{v},
\end{eqnarray}
where $\Vec{v}$ in Eq.~(\ref{dnI}) is given in Eq.~(\ref{vI}). Besides, from Eq.~(\ref{wK}),
\begin{equation} \label{dnu}
\dot\theta^*=\frac{G}{r^2}.
\end{equation}
Finally, the variation of $G$ is obtained by replacing Eq.~(\ref{vI}) into Eq.~(\ref{Gp}), and the variation of $r$ by replacing Eq.~(\ref{uI}) into Eq.~(\ref{ddr}).
\par

Now, the differential system given by Eqs.~(\ref{Gp}), (\ref{ddr}), and (\ref{duI})--(\ref{dnu}) is of dimension thirteen with the same six constraints in Eq.~(\ref{constraints}) plus the additional constraint given by $\dot{\Vec{u}}^*\cos\theta^*+\dot{\Vec{v}}^*\sin\theta^*=\Vec{0}$ immediately obtained from Eqs.~(\ref{duI}) and (\ref{dvI}), which is in fact a scalar constraint because both $\dot{\Vec{u}}^*$ and $\dot{\Vec{v}}^*$ have the direction of $\Vec{n}$.
\par

In order to evaluate the components of the disturbing force in the directions of the ideal frame, at each integration step $t=t_i$ the direction vectors of the ideal frame must be projected onto the space frame,
\begin{equation}
(\Vec{u}_\mathcal{S}^*(t_i),\Vec{v}_\mathcal{S}^*(t_i),\Vec{n}_\mathcal{S}(t_i))
=\mathcal{M}(t_0)\mathcal{N}(t_i)
\end{equation}
where $\mathcal{M}$ is given in Eq.~(\ref{SO}), and the time dependent rotation matrix
\begin{equation} \label{Qm}
\mathcal{N}\equiv(\Vec{u}_\mathcal{D}^*(t),\Vec{v}_\mathcal{D}^*(t),\Vec{n}_\mathcal{D}(t))
\end{equation}
provides the components of the ideal frame in the departure frame. 
\par

\subsection{Euler parameters} \label{ss:Eulerp}

The dimension of the differential system is reduced by elimination of constraints. It can be done representing the rotation of the ideal frame by the Eulerian angles, namely $\Omega^*$, for the argument of the node, $I^*$, for inclination, and $\sigma$, for the longitude of $\Vec{u}^*$ reckoned from the node. The variation of the Eulerian angles is obtained from the usual relations \cite[p.~4]{Leimanis1965}
\begin{equation}\label{wb}
\Vec\omega^*=\left(\begin{array}{c} 
\dot\Omega^*\sin{I}^*\sin\sigma+\dot{I}^*\cos\sigma \\ 
\dot\Omega^*\sin{I}^*\cos\sigma-\dot{I}^*\sin\sigma \\ 
\dot\Omega^*\cos{I}^*+\dot\sigma
\end{array}\right).
\end{equation}
Hence, in view of Eq.~(\ref{wIss}),
\begin{eqnarray} \label{hp}
\dot\Omega^* &=& \frac{r}{G}(\Vec{P}\cdot\Vec{n})\frac{\sin(\theta^*+\sigma)}{\sin{I}}, \\ \label{ip}
\dot{I}^* &=& \frac{r}{G}(\Vec{P}\cdot\Vec{n})\cos(\theta^*+\sigma), \\ \label{sp}
\dot\sigma &=& -\dot\Omega^*\cos{I},
\end{eqnarray}
which admit the non-holonomic constraint
\begin{equation} \label{bilinearEulerAngles}
\dot{I}^*\sin(\theta^*+\sigma)-\dot\Omega^*\sin{I}^*\cos(\theta^*+\sigma)=0.
\end{equation}
\par

However, this differential system is singular for $I^*=0$, which corresponds to the initial conditions in the departure frame: 
$\Omega^*(t_0)=I^*(t_0)=\sigma(t_0)=0$. Among the different non-singular variables that can be chosen to avoid the singularity, a notable case is given by the set of Euler redundant parameters
\begin{eqnarray} \label{q1}
\lambda_1 &=& \sin\mbox{$\frac{1}{2}$}I^*\cos\mbox{$\frac{1}{2}$}(\Omega^*-\sigma), \\ \label{q2}
\lambda_2 &=& \sin\mbox{$\frac{1}{2}$}I^*\sin\mbox{$\frac{1}{2}$}(\Omega^*-\sigma), \\ \label{q3}
\lambda_3 &=& \cos\mbox{$\frac{1}{2}$}I^*\sin\mbox{$\frac{1}{2}$}(\Omega^*+\sigma), \\ \label{q0}
\lambda_4 &=& \cos\mbox{$\frac{1}{2}$}I^*\cos\mbox{$\frac{1}{2}$}(\Omega^*+\sigma), 
\end{eqnarray}
with the geometric constraint
\begin{equation} \label{cq1}
\lambda_1^2+\lambda_2^2+\lambda_3^2+\lambda_4^2=1.
\end{equation}
Then, as easily checked, the rotation matrix $\mathcal{N}=R_3(-\Omega^*)\,R_1(-I^*)\,R_3(-\sigma)$ in Eq.~(\ref{Qm}), in which $R_1$ and $R_3$ are the usual rotation matrices about the axis $x$ and $z$, respectively, is written
\begin{equation} \label{rotaQ}
\mathcal{N}=\left(
\begin{array}{ccc}
 1-2(\lambda_2^2+\lambda_3^2) & 2 \left(\lambda_1 \lambda_2-\lambda_4 \lambda_3\right) & 2 \left(\lambda_1 \lambda_3+\lambda_4 \lambda_2\right) \\
 2 \left(\lambda_1 \lambda_2+\lambda_4 \lambda_3\right) & 1-2(\lambda_1^2+\lambda_3^2) & 2 \left(\lambda_2 \lambda_3-\lambda_4 \lambda_1\right) \\
 2 \left(\lambda_1 \lambda_3-\lambda_4 \lambda_2\right) & 2 \left(\lambda_2 \lambda_3+\lambda_4 \lambda_1\right) & 1-2(\lambda_1^2+\lambda_2^2)
\end{array}
\right).
\end{equation}
\par

Differentiation of Eqs.~(\ref{q1})--(\ref{q0}) with the concomitant use of Eqs.~(\ref{hp})--(\ref{sp}) yields
\begin{eqnarray} \label{qp1}
\dot\lambda_1 &=& \frac{r}{2G}(\Vec{P}\cdot\Vec{n})(\lambda_4\cos\theta^*-\lambda_3\sin\theta^*), \\ \label{qp2}
\dot\lambda_2 &=& \frac{r}{2G}(\Vec{P}\cdot\Vec{n})(\lambda_4\sin\theta^*+\lambda_3\cos\theta^*), \\ \label{qp3}
\dot\lambda_3 &=& \frac{r}{2G}(\Vec{P}\cdot\Vec{n})(\lambda_1\sin\theta^*-\lambda_2\cos\theta^*), \\ \label{qp0}
\dot\lambda_4 &=& \frac{r}{2G}(\Vec{P}\cdot\Vec{n})(-\lambda_1\cos\theta^*-\lambda_2\sin\theta^*),
\end{eqnarray}
with the non-holonomic constraint
\begin{equation} \label{cq2}
\lambda_1\dot\lambda_2-\lambda_2\dot\lambda_1+\lambda_3\dot\lambda_4-\lambda_4\dot\lambda_3=0,
\end{equation}
which is the analog of Eq.~(\ref{bilinearEulerAngles}). Equations (\ref{qp1})--(\ref{qp0}) are conveniently used in replacement of Eqs.~(\ref{duI})--(\ref{dnI}).
\par

\subsection{A modified set of Eulerian parameters} \label{ss:modEulerp}

Alternatively to the Eulerian parameter in Eqs.~(\ref{q1})--(\ref{q0}), the modified set
\begin{equation} \label{gi}
g_i=\sqrt{G}\lambda_i, \qquad i=1,\dots4,
\end{equation}
can be used. The geometric constraint in Eq.~(\ref{cq1}) no longer applies, being replaced by
\begin{equation} \label{Ggg}
G=g_1^2+g_2^2+g_3^2+g_4^2,
\end{equation}
which shows that the modulus of the angular momentum is no longer a variable but a derived quantity.

Differentiation of Eq.~(\ref{gi}), using Eqs.~(\ref{qp1})--(\ref{qp0}) and replacing $\dot{G}$ by the right side of Eq.~(\ref{Gp}), yields
\begin{eqnarray} \label{gp1}
\dot{g}_1 &=& \frac{r}{2G}\left[(\Vec{P}\cdot\Vec{v})g_1+(\Vec{P}\cdot\Vec{n})(g_4\cos\theta^*-g_3\sin\theta^*)\right], \\ \label{gp2}
\dot{g}_2 &=& \frac{r}{2G}\left[(\Vec{P}\cdot\Vec{v})g_2+(\Vec{P}\cdot\Vec{n})(g_4\sin\theta^*+g_3\cos\theta^*)\right], \\ \label{gp3}
\dot{g}_3 &=& \frac{r}{2G}\left[(\Vec{P}\cdot\Vec{v})g_3+(\Vec{P}\cdot\Vec{n})(g_1\sin\theta^*-g_2\cos\theta^*)\right], \\ \label{gp0}
\dot{g}_4 &=& \frac{r}{2G}\left[(\Vec{P}\cdot\Vec{v})g_4-(\Vec{P}\cdot\Vec{n})(g_1\cos\theta^*+g_2\sin\theta^*)\right].
\end{eqnarray}
\par

The differential system is now of dimension 7, comprising Eqs.~(\ref{ddr}), (\ref{dnu}), and (\ref{gp1})--(\ref{gp0}), and it is easy to check that Eq.~(\ref{cq2}) still applies when replacing $\lambda_i$ by corresponding $g_i$, viz.
\begin{equation} \label{cggp}
g_1\dot{g}_2-g_2\dot{g}_1+g_3\dot{g}_4-g_4\dot{g}_3=0.
\end{equation}

\subsection{The ellipse in the orbital plane}

The position of the osculating ellipse in the orbital plane is described by the attitude of the apsidal frame with respect to the ideal frame. Because $\Vec{e}$ is undefined for circular orbits, it is convenient to use Deprit's (\citeyear{Deprit1975}) \emph{ideal elements}
\[
C^*=\frac{G}{p}\Vec{e}\cdot\Vec{u}^* =\frac{G}{p}e\cos\gamma, \qquad 
S^*=\frac{G}{p}\Vec{e}\cdot\Vec{v}^*=\frac{G}{p}e\sin\gamma,
\]
where the angle $\gamma$ is reckoned from $\Vec{u}^*$ to $\Vec{e}$ counterclockwise. Hence 
\[
C\Vec{u}^*+S\Vec{v}^*=\frac{G}{p}\Vec{e},
\]
and standard differentiation in the ideal frame, in which $\partial\Vec{e}/\partial{t}=\dot{\Vec{e}}-\Vec\omega^*\times\Vec{e}$, from Eq.~(\ref{darboux}) with $\Vec\omega^*$ is given by Eq.~(\ref{wIss}), yields
\begin{equation} \label{CuSup}
\dot{C}^*\Vec{u}^*+\dot{S}^*\Vec{v}^* =\left(1+\frac{r}{p}\right)(\Vec{P}\cdot\Vec{v})\Vec{u}-(\Vec{P}\cdot\Vec{u})\Vec{v},
\end{equation}
where $\Vec{u}$ and $\Vec{v}$ are given in Eqs.~(\ref{uI}) and (\ref{vI}), respectively. Multiplication of both sides of Eq.~(\ref{CuSup}) by $\Vec{u}^*$ and $\Vec{v}^*$, respectively, yields
\begin{eqnarray} \label{dC}
\dot{C}^* &=& \left(1+\frac{r}{p}\right)(\Vec{P}\cdot\Vec{v})\cos\theta^* + (\Vec{P}\cdot\Vec{u})\sin\theta^*, 
\\ \label{dS}
\dot{S}^* &=& \left(1+\frac{r}{p}\right)(\Vec{P}\cdot\Vec{v})\sin\theta^* - (\Vec{P}\cdot\Vec{u})\cos\theta^*.
\end{eqnarray}
\par

In view of the relations in the Keplerian ellipse
\[
r=\frac{p}{1+e\cos(\theta^*-\gamma)}, \qquad \dot{r}=\frac{G}{p}e\sin(\theta^*-\gamma),
\]
one easily gets
\begin{eqnarray} \label{qCS}
\frac{G}{r} &=& C^*\cos\theta^*+S^*\sin\theta^*+\frac{G}{p},
\\ \label{qqCS}
\dot{r} &=& C^*\sin\theta^*-S^*\cos\theta^*,
\end{eqnarray}
and hence
\begin{eqnarray} \label{CrR}
C^* &=& \left(\frac{G}{r}-\frac{G}{p}\right)\cos\theta^*+\dot{r}\sin\theta^*, \\ \label{SrR}
S^* &=& \left(\frac{G}{r}-\frac{G}{p}\right)\sin\theta^*-\dot{r}\cos\theta^*.
\end{eqnarray}
Therefore $C^*$ and $S^*$ are functions of $r$ and $\dot{r}$, and Eqs.~(\ref{dC})--(\ref{dS}) can be used to replace Eq.~(\ref{ddr}) so that the differential system to integrate depends on either $(\lambda_1,\lambda_2,\lambda_3,\lambda_4,G,C^*,S^*)$ or  $(g_1,g_2,g_3,g_4,C^*,S^*)$, all of which vary slowly, and $\theta^*$. 
\par

At each step of the integration, the Cartesian coordinates are computed from $\Vec{x}=r\Vec{u}$ and $\Vec{X}=\dot{r}\Vec{u}+r\dot{\Vec{u}}$, with $\dot{\Vec{u}}$ given in Eq.~(\ref{vup}). That is,
\begin{eqnarray} \label{xP}
\Vec{x} &=& r\,\mathcal{M}(t_0)\mathcal{N}(t)\Vec{u}_\mathcal{I}, \\ \label{xQ}
\Vec{X} &=& \frac{\dot{r}}{r}\Vec{x} + \frac{G}{r}\mathcal{M}(t_0)\mathcal{N}(t)\Vec{v}_\mathcal{I}, 
\end{eqnarray}
where, from Eqs.~(\ref{uI}) and (\ref{vI}),
\begin{equation} \label{uIvI}
\Vec{u}_\mathcal{I}=\left(\begin{array}{c} \cos\theta^* \\ \sin\theta^* \\ 0 \end{array}\right), \qquad
\Vec{v}_\mathcal{I}=\left(\begin{array}{c} -\sin\theta^* \\ \cos\theta^* \\ 0 \end{array}\right).
\end{equation}
\par

\subsection{Regularization}

The explicit appearance of $r$ in denominators of the Keplerian terms of Eq.~(\ref{ddr}), as well in Eq.~(\ref{dnu}), may harm the numerical integration of highly elliptic orbits close to the periapsis. These undesired denominators can be removed by regularizing the differential system making the change $q=1/r$ and using the angle $\theta^*$ as the new time defined by Eq.~(\ref{dnu}), from which
\[
\frac{\mathrm{d}}{\mathrm{d}\theta^*} =\frac{r^2}{G}\frac{\mathrm{d}}{\mathrm{d}t} =\frac{1}{q^2G}\frac{\mathrm{d}}{\mathrm{d}t}.
\]
Hence, standard operations transform Eqs.~(\ref{gp1})--(\ref{gp0}) into
\begin{eqnarray} \label{qp1r}
\frac{\mathrm{d}g_1}{\mathrm{d}\theta^*} &=& \frac{1}{2}\left[(\Vec{P}^*\cdot\Vec{v})g_1+(\Vec{P}^*\cdot\Vec{n})(g_4\cos\theta^*-g_3\sin\theta^*)\right], \\ 
\label{qp2r}
\frac{\mathrm{d}g_2}{\mathrm{d}\theta^*} &=& \frac{1}{2}\left[(\Vec{P}^*\cdot\Vec{v})g_2+(\Vec{P}^*\cdot\Vec{n})(g_4\sin\theta^*+g_3\cos\theta^*)\right], \\ \label{qp3r}
\frac{\mathrm{d}g_3}{\mathrm{d}\theta^*} &=& \frac{1}{2}\left[(\Vec{P}^*\cdot\Vec{v})g_3+(\Vec{P}^*\cdot\Vec{n})(g_1\sin\theta^*-g_2\cos\theta^*)\right], \\ \label{qp0r}
\frac{\mathrm{d}g_4}{\mathrm{d}\theta^*} &=& \frac{1}{2}\left[(\Vec{P}^*\cdot\Vec{v})g_4-(\Vec{P}^*\cdot\Vec{n})(g_1\cos\theta^*+g_2\sin\theta^*)\right], \\ 
\label{rpr}
\frac{\mathrm{d}q}{\mathrm{d}\theta^*}   &=& Q, \\ \label{ddrIr}
\frac{\mathrm{d}Q}{\mathrm{d}\theta^*}   &=& \frac{1}{p}-q\left[1+(\Vec{P}^*\cdot\Vec{u})\right] -Q(\Vec{P}^*\cdot\Vec{v}), \\ \label{tr}
\frac{\mathrm{d}t}{\mathrm{d}\theta^*}   &=& \frac{1}{q^2G},
\end{eqnarray}
where $Q=-\dot{r}/G$, and the non-dimensional force $\Vec{P}^*\equiv\Vec{P}/(q^3G^2)$ is used as abbreviation. The constraint in Eq.~(\ref{cggp}) still applies by trivially changing time differentiation by differentiation with respect to $\theta^*$. 
\par

The modified Eulerian parameters are elements which vary slowly. Besides,  when the perturbations vanish the differential system reduces to the integration of Eqs.~(\ref{rpr}) and (\ref{ddrIr}), which become linear resulting in the harmonic oscillations $\mathrm{d}^2q/\mathrm{d}{\theta^{*2}}=-q+\mathcal{G}M/G^2$.
\par

Algorithm \ref{al:qQ} summarizes the necessary operations for implementing the integration of the regularized Eqs.~(\ref{qp1r})--(\ref{tr}), and shows the simplicity of this approach. In particular, like with Cowell's formulation, the perturbation can be evaluated directly in the inertial frame . Besides, it is customary to use internal units of length and time, which are commonly chosen as $\mathrm{UL}=\mathcal{G}M/(-\Vec{X}\cdot\Vec{X}+2\mathcal{G}M/r)$, and $\mathrm{UT}=\mathrm{UL}\sqrt{\mathrm{UL}/\mathcal{G}M}$.
\par

\begin{algorithm}[htbp]
\caption{\label{al:qQ} Integration of Eqs.~(\ref{qp1r})--(\ref{tr})}
\begin{algorithmic}[1] 
\State \textbf{Inputs}: Initial epoch $t_0$, final epoch $T$, evaluation interval $\Delta\theta^*$;  $\Vec{x}_0$, $\Vec{X}_0$
\State call \Call{initial conditions}{$t_0,\Vec{x}_0,\Vec{X}_0,g_1,g_2,g_3,g_4,q,Q$}
\While{$t\le{T}$ }
\State call \Call{force model}{$\theta^*,t,g_1,g_2,g_3,g_4,q,Q,\Vec{P}^*\cdot\Vec{u}_\mathcal{S},\Vec{P}^*\cdot\Vec{v}_\mathcal{S},\Vec{P}^*\cdot\Vec{n}_\mathcal{S}$}
\State numerically integrate Eqs.~(\ref{qp1r})--(\ref{tr}) with $1/p=\mathcal{G}M/G^2$,
\If{output required}
\State make $r=1/q$, $\dot{r}=-QG$, where $G$ is computed from Eq.~(\ref{Ggg}), 
\State evaluate Eqs.~(\ref{xP}) and (\ref{xQ}) using Eq.~(\ref{uIvI}); save $t$, $\Vec{x}$, $\Vec{X}$
\EndIf
\State $\theta^*\leftarrow\theta^*+\Delta\theta^*$
\EndWhile
\Statex
\Procedure{initial conditions}{$t,\Vec{x},\Vec{X},g_1,g_2,g_3,g_4,q,Q$}
\State compute $\Vec{G}$ from Eq.~(\ref{xxX}), and $\Vec{u}$, $\Vec{n}$, $\Vec{v}$ from Eq.~(\ref{orbitalframe}) with $G=\|\Vec{G}\|$, $r=\|\Vec{x}\|$
\State compute the constant matrix $\mathcal{M}(t_0)$ using Eq.~(\ref{SO})
\State make $\theta^*(t)=0$, $g_1=g_2=g_3=0$, $g_4=\sqrt{G}$; compute $\dot{r}=\Vec{X}\cdot\Vec{x}/r$ 
\State make $q=1/r$, $Q=-\dot{r}/G$
\EndProcedure
\Statex
\Procedure{force model}{$\theta^*,t,g_1,g_2,g_3,g_4,q,Q,\Vec{P}^*\cdot\Vec{u}_\mathcal{S},\Vec{P}^*\cdot\Vec{v}_\mathcal{S},\Vec{P}^*\cdot\Vec{n}_\mathcal{S}$}
\State compute $G$ from Eq.~(\ref{Ggg}), make $\lambda_i=g_i/\sqrt{G}$, compute $\mathcal{N}(t)$ from Eq.~(\ref{rotaQ})
\State compute $\Vec{u}_\mathcal{I}$, $\Vec{v}_\mathcal{I}$ from Eq.~(\ref{uIvI}); make $\Vec{n}_\mathcal{I}=(0,0,1)$
\State compute $(\Vec{u}_\mathcal{S},\Vec{v}_\mathcal{S},\Vec{n}_\mathcal{S})=\mathcal{M}(t_0)\,\mathcal{N}(t)\,(\Vec{u}_\mathcal{I},\Vec{v}_\mathcal{I},\Vec{n}_\mathcal{I})$
\State evaluate $\Vec{P}_\mathcal{S}$ in the space frame; it depends on the problem at hand
\State scale $\Vec{P}^*\equiv(r^3/G^2)\Vec{P}_\mathcal{S}$ and compute $(\Vec{P}^*\cdot\Vec{u}_\mathcal{S})$, $(\Vec{P}^*\cdot\Vec{v}_\mathcal{S})$, $(\Vec{P}^*\cdot\Vec{n}_\mathcal{S})$
\EndProcedure
\end{algorithmic}
\end{algorithm}

On the other hand, Eqs.~(\ref{rpr}) and (\ref{ddrIr}) can be replaced by
\begin{eqnarray} \label{dCnu}
\frac{\mathrm{d}C^*}{\mathrm{d}\theta^*} &=& 
\left(\frac{G}{r}+\frac{G}{p}\right)(\Vec{P}^*\cdot\Vec{v})\cos\theta^*+\frac{G}{r}(\Vec{P}^*\cdot\Vec{u})\sin\theta^*, \\ \label{dSnu}
\frac{\mathrm{d}S^*}{\mathrm{d}\theta^*} &=& 
\left(\frac{G}{r}+\frac{G}{p}\right)(\Vec{P}^*\cdot\Vec{v})\sin\theta^*-\frac{G}{r}(\Vec{P}^*\cdot\Vec{u})\cos\theta^*,
\end{eqnarray}
which are trivially obtained from Eqs.~(\ref{dC}) and (\ref{dS}), respectively, where $G/r$ is given in Eq.~(\ref{qCS}) and $G/p=\mathcal{G}M/G$. The variations introduced by this second approach are illustrated with Algorithm \ref{al:CS}. The procedure for evaluating the projections of the disturbing force in the ideal frame is the same as in Algorithm \ref{al:qQ}, except for the use of Eqs.~(\ref{qCS}) and (\ref{qqCS}), and is not presented.

\par

\begin{algorithm}[htbp]
\caption{\label{al:CS} Integration of Deprit's ideal elements}
\begin{algorithmic}[1] 
\State \textbf{Inputs}: Initial epoch $t_0$, final epoch $T$, evaluation interval $\Delta\theta^*$;  $\Vec{x}_0$, $\Vec{X}_0$
\State call \Call{initial conditions}{$t_0,\Vec{x}_0,\Vec{X}_0,g_1,g_2,g_3,g_4,C^*,S^*$}
\While{$t\le{T}$ }
\State call \Call{force model}{$\theta^*,t,g_1,g_2,g_3,g_4,C^*,S^*,\Vec{P}^*\cdot\Vec{u}_\mathcal{S},\Vec{P}^*\cdot\Vec{v}_\mathcal{S},\Vec{P}^*\cdot\Vec{n}_\mathcal{S}$}
\State compute $G/r$ from Eq.~(\ref{qCS}), make $G/p=\mathcal{G}M/G$
\State numerically integrate Eqs.~(\ref{qp1r})--(\ref{qp0r}), and (\ref{tr})--(\ref{dSnu})
\If{output required}
\State compute $G$ from Eq.~(\ref{Ggg}), $r$ from Eq.~(\ref{qCS}) and $\dot{r}$ from Eq.~(\ref{qqCS})
\State evaluate Eqs.~(\ref{xP}) and (\ref{xQ}), save $t_i$, $\Vec{x}_i$, $\Vec{X}_i$
\EndIf
\State $\theta^*\leftarrow\theta^*+\Delta\theta^*$
\EndWhile
\Statex
\Procedure{initial conditions}{$t,\Vec{x},\Vec{X},g_1,g_2,g_3,g_4,C^*,S^*$}
\State same as lines 13--15 of Algorithm \ref{al:qQ}
\State make $G/p=\mathcal{G}M/G$, evaluate $C^*$ and $S^*$ from Eq.~(\ref{CrR}) and (\ref{SrR})
\EndProcedure
\end{algorithmic}
\end{algorithm}

\section{Performance evaluation} \label{s:performance}

The efficiency of the new formulation of the equations of motion has been tested by comparing accuracy and runtime with respect to the performance of the classical formulation of dimension 8, which integrates Eqs.~(\ref{Gp}) and (\ref{qp1})--(\ref{qp0}) instead of Eqs.~(\ref{gp1})--(\ref{gp0}). The advantages of using an integration method in associattion with a particular formulation may depend on the orbital scenario to which the formulation is applied \citep[see chap.~6 of][]{RoaThesis2016}, and the discussion of the more efficient integration method for the ideal frame formulation is not approached here, where both differential systems have been integrated numerically with the reliable and widespread \href{http://www.unige.ch/~hairer/prog/nonstiff/dop853.f}{DOP853} free code described in \citep{HairerNorsettWanner2008}.
\par

Following tradition, the tests were based on the numerical examples in \cite[p.~118 and ff.]{StiefelScheifele1971}. Namely, the forces model considers the non-centralities of the Geopotential limited to the contribution of the second zonal harmonic, as well as the moon attraction in the simplifying assumption that the moon moves in a circular orbit about the earth. The ``true'', reference orbit was borrowed from \citep{UrrutxuaSanjurjoRivoPelaez2016}, who used extended precision to retain the common figures of the different solutions obtained with a variety of high precision integrators.
\par

The 7 dimensional formulation was always found to improve performance with respect to the 8 dimensional case in terms of computing time by about 3\%, while retaining the same accuracy. On the other hand, the integration of the inverse of the distance variant shows a little bit faster than the integration of the ideal elements $C^*$ and $S^*$, but at the expense of a slightly loss of accuracy.
\par

Finally, it worths mentioning that the performance of the non-regularized version based on the integration of the ideal elements $C^*$ and $S^*$ is only slightly worse than the regularized version in terms of accuracy, yet, as expected, it is penalized by generally doubling runtime. On the contrary, in addition to the important increase of runtime, the accuracy of the integration notably deteriorates when integrating the distance and radial velocity without regularization.

\section{Conclusions}
Simple modifications in the definition of the Eulerian parameters lead to a new ideal frame formalism with a minimum dimension rising from 6 to 7 dimensions. The new formulation takes advantage of standard regularization, is conceptually very simple, and enjoys slightly better performance than its traditional counterpart in dimension 8. The higher accuracy is always obtained when the motion in the orbital plane is materialized by the integration of Deprit's ideal elements, which evolve slowly either in the physical or regularized time scales. On the other hand, the integration of the inverse of the radial distance is simpler and faster than that of the ideal elements, although it performs slightly worse in terms of accuracy.

\section*{Acknowledgemnts}

This work is partially supported by the Ministry of Economic Affairs and Competitiveness of Spain, under grants ESP2013-41634-P and ESP2014-57071-R. The author thanks S.~Ferrer, University of Murcia, for his comments on a preliminary manuscript.

\end{document}